\documentclass[12pt]{article}
\usepackage{amsfonts,amssymb}

\def\C{\mathbb{C}}
\def\Z{\mathbb{Z}}
\def\N{\mathbb{N}}

\def\bq{ \begin{equation} }
\def\eq{ \end{equation} }
\def\ben{ \begin{eqnarray} }
\def\en{ \end{eqnarray} }

\def\frac#1#2{{#1\over #2}}

\def\on#1#2{\mathop{\vbox{\ialign{##\crcr\noalign{\kern2pt}
$\scriptstyle{#2}$\crcr\noalign{\kern2pt\nointerlineskip}
\kern-2pt$\hfil\displaystyle{#1}\hfil$\crcr}}}\limits}

\title{{\bf Compatible quadratic Poisson brackets 
related to a family of elliptic curves}}
\author{Alexander Odesskii and Thomas Wolf\\ 
Department of Mathematics,
Brock University\\ 500 Glenridge Avenue, St.Catharines, 
Ontario, Canada L2S 3A1\\
email: aodesski@brocku.ca, twolf@brocku.ca}

\begin{document}
\maketitle
\begin{abstract}
We construct nine pairwise compatible quadratic Poisson structures
such that a generic linear combination of them is associated with an
elliptic algebra in $n$ generators. Explicit formulas for Casimir elements of this elliptic Poisson structure are obtained. 
\vspace{1cm}  
\end{abstract}

%

%
%
%
%
\newpage
\tableofcontents
\newpage

\section{Introduction} \label{intro}
Two Poisson brackets $\{\cdot,\cdot\}_0$ and $\{\cdot,\cdot\}_1$ defined on the same
finite dimensional vector space are said to be compatible if 
\begin{equation} \label{pensil}
\{\cdot,\cdot\}_{u}=\{\cdot,\cdot\}_0+u \{\cdot,\cdot\}_1
\end{equation}
is a Poisson bracket for any constant $u$. Note that if
$\{\cdot,\cdot\}_{u_1,...,u_k}=\{\cdot,\cdot\}_0+u_1
\{\cdot,\cdot\}_1+...+u_k \{\cdot,\cdot\}_k$ is a Poisson bracket for
arbitrary $u_1,...,u_k$, then all brackets
$\{\cdot,\cdot\}_0,...,\{\cdot,\cdot\}_k$ are Poisson and pairwise
compatible. Compatible Poisson structures play an important role in
the theory of integrable systems \cite{magri,bols} and in differential
geometry \cite{gelzah,gelzah1}. A lot of examples of compatible
Poisson structures are known \cite{bols}. Most of these are linear in
certain coordinates. However, quadratic Poisson structures are also
interesting. While the theory of linear Poisson structures is
well-understood and possesses a classification theory\footnote{The
  theory of linear Poisson structures coincides with the theory of Lie
  algebras.}, the theory of quadratic Poisson algebras is more
complicated. If the dimension of a linear space is larger than four,
then no classification results for quadratic Poisson structures on
this space are available.  All known examples can be divided into two
classes: rational and elliptic.  In the elliptic case structure
constants of a Poisson bracket are modular functions of a parameter
$\tau$, a modular parameter of an elliptic curve. This elliptic curve
appears naturally as a symplectic leaf of this elliptic Poisson
structure \cite{odes}.

Let $Q_n(\tau,\eta)$ be an associative algebra defined by $n$
generators $\{x_i;i\in\Z/n\Z\}$ and quadratic relations \cite{odes}
$$\sum_{r\in\Z/n\Z}\frac{1}{\theta_{j-i-r}(-\eta,\tau)\theta_r(\eta,\tau)}x_{j-r}x_{i+r}=0,$$
for all $i\ne j\in\Z/n\Z$. Here $\theta_i(z,\tau)$ are
$\theta$-functions with characteristics (see Appendix). It is known
that for generic $\eta$ the algebra $Q_n(\tau,\eta)$ has the same size
of graded components as the polynomial ring
$\C[x_1,...,x_n]$. Moreover, if $\eta=0$, then $Q_n(\tau,\eta)$ is
isomorphic to $\C[x_1,...,x_n]$. Therefore, for any fixed $\tau$ we
have a flat deformation of a polynomial ring. Let $q_n(\tau)$ be the
corresponding Poisson algebra. Symplectic leaves of this Poisson
structure are known \cite{odes}. In particular, the center of this
Poisson algebra is generated by one homogeneous polynomial of degree
$n$ if $n$ is odd and by two homogeneous polynomials of degree
$\frac{n}{2}$ if $n$ is even.

One can pose the following problems: 
\begin{enumerate}
\item Do there exist Poisson structures compatible with the one in
$q_n(\tau)$?
\item Construct a maximal number of Poisson structures pairwise
compatible and compatible with the one in $q_n(\tau)$.
\end{enumerate}
It is easy to study these questions in the cases $n=3,~4$.

Let $n=3$. The Poisson bracket in $q_3(\tau)$ can be written as 
$$\{x_{\sigma_1},x_{\sigma_2}\}=\frac{\partial P}{\partial
  x_{\sigma_3}}$$ where $P$ is a certain homogeneous cubic polynomial
in $x_1,x_2,x_3$ and $\sigma$ is an arbitrary even
permutation. Moreover, this formula defines a Poisson bracket for an
arbitrary polynomial $P$ and all these brackets are pairwise
compatible. In particular, there exist 10 linearly independent
quadratic Poisson brackets because there are 10 linearly independent
homogeneous cubic polynomials in 3 variables.

Let $n=4$. The Poisson bracket in $q_4(\tau)$ can be written as
$$\{x_{\sigma_1},x_{\sigma_2}\}=\det\left(\begin{array}{cc}
\frac{\partial P}{\partial x_{\sigma_3}}&
\frac{\partial P}{\partial x_{\sigma_4}} \vspace{6pt}\\
\frac{\partial R}{\partial x_{\sigma_3}}&
\frac{\partial R}{\partial x_{\sigma_4}}
\end{array}\right)$$
where $P,~R$ are certain homogeneous quadratic polynomials in
$x_1,...,x_4$ and $\sigma$ is an arbitrary even permutation. Moreover,
this formula defines a Poisson bracket for arbitrary polynomials $P$
and $R$. If we fix $R$ and vary $P$, we obtain an infinite family of
pairwise compatible Poisson brackets. In particular, there exist 9
pairwise compatible quadratic Poisson brackets.  Indeed, there are 10
quadratic polynomials in 4 variables and $P$ should not be
proportional to $R$.

If $n>4$, then the situation is more complicated because the similar
construction for $q_n(\tau)$, $n>4$ does not exists. In this paper we
construct nine pairwise compatible quadratic Poisson brackets
\footnote{Three of these where constructed in \cite{od1}.}  for
arbitrary $n$. A generic linear combination of these Poisson brackets
is isomorphic to $q_n(\tau)$ where $\tau$ depends on coefficients in
this linear combination. Moreover, we think that this family of
Poisson brackets is maximal. We have checked, that for $n=5,6,...,40$
there are no quadratic Poisson brackets that are compatible with all
our nine Poisson brackets and are linearly independent of them. For
these values of $n$ there are no Poisson brackets compatible with all ours that are constant, linear, cubic and quartic.

Let us describe the contents of the paper. In section \ref{algebraic}
we construct nine compatible quadratic Poisson structures on a certain
$n$-dimensional linear space ${\cal F}_n$. This construction is slightly different for even and odd $n$. It is 
summarized in Remarks 1, 1$^{\prime}$ and 2, 2$^{\prime}$ as an algorithm for the computation of
Poisson brackets between $x_i$ and $x_j$ for $i,j=0,2,3,...,n$. In
section \ref{functional} we explain the functional version of the same
construction. In sections \ref{cas} we describe symplectic leaves and Casimir elements of our Poisson algebras (see also \cite{odes, od1}). In the Conclusion we outline several open problems. In the
Appendix we collect some notations and standard facts about elliptic
and $\theta$-functions (see \cite{ell,odes} for details).

\section{Algebraic construction of nine compatible quadratic 
         Poisson brackets} \label{algebraic}

\subsection{Notations}

Our construction of Poisson brackets is slightly different for even and odd $n$. We will use index $ev$ (resp. $od$) for objects related to even (resp. odd) $n$.

Let 
\begin{equation}\label{pol}
P_{ev}(t)=a_0+a_1t+a_2t^2+a_3t^3+a_4t^4,~~P_{od}(t)=a_0+a_1t+a_2t^2+a_3t^3,~~Q(t)=b_0+b_1t+b_2t^2 
\end{equation}
 be arbitrary polynomials of degree not larger then four, three and two correspondingly. Let ${\cal F}_{ev}$ be a commutative 
associative algebra defined by generators $f,~g$ and the relation 
\begin{equation}\label{relev}
g^2=P_{ev}(f)+Q(f)g. 
\end{equation}
Let ${\cal F}_{od}$ be a commutative 
associative algebra defined by generators $f,~g$ and the relation 
\begin{equation}\label{relod}
(f+c)g^2=P_{od}(f)+Q(f)g. 
\end{equation}
Here $a_0,...,a_4,b_0,b_1,b_2,c$ are constants. Let $D$ be a derivation of ${\cal F}_{ev}$ and ${\cal F}_{od}$ defined on its generators by
\begin{equation}\label{derev}
D(f)=2g-Q(f),~~~D(g)=P_{ev}^{\prime}(f)+Q^{\prime}(f)g
\end{equation}
for ${\cal F}_{ev}$ and by
\begin{equation}\label{derod}
D(f)=2(f+c)g-Q(f),~~~D(g)=P_{od}^{\prime}(f)+Q^{\prime}(f)g-g^2
\end{equation}
for ${\cal F}_{od}$.

Let $\cal F$ be either ${\cal F}_{ev}$ or ${\cal F}_{od}$.
It is clear that $\cal F\otimes\cal F$ is generated by $f_1=f\otimes 1,~f_2=1\otimes f,~g_1=g\otimes 1,~g_2=1\otimes g$ as an associative algebra. For an arbitrary element $h\in\cal F$ 
we will use the notations $h_1=h\otimes 1,~h_2=1\otimes h$ for the corresponding elements in $\cal F\otimes\cal F$. Let $\lambda_{ev}\in Frac({\cal F}_{ev}\otimes{\cal F}_{ev})$ be an element of a field of fractions of 
${\cal F}_{ev}\otimes{\cal F}_{ev}$ 
defined by
\begin{equation}\label{homev}
(f_1-f_2)\lambda_{ev}=g_1+g_2-\frac{1}{2}Q(f_1)-\frac{1}{2}Q(f_2),
\end{equation}
or by
$$(g_1-g_2)\lambda_{ev}=\frac{P_{ev}(f_1)-P_{ev}(f_2)}{f_1-f_2}+\frac{Q(f_1)-Q(f_2)}{2(f_1-f_2)}(g_1+g_2).$$
These definitions are equivalent by virtue of (\ref{relev}).
Let $\lambda_{od}\in Frac({\cal }F_{od}\otimes{\cal F}_{od})$ be an element of a field of fractions of ${\cal F}_{od}\otimes{\cal F}_{od}$ 
defined by
\begin{equation}\label{homod}
(f_1-f_2)\lambda_{od}=(f_1+c)g_1+(f_2+c)g_2-\frac{1}{2}Q(f_1)-\frac{1}{2}Q(f_2),
\end{equation}
or by
$$(g_1-g_2)\lambda_{od}=\frac{P_{od}(f_1)-P_{od}(f_2)}{f_1-f_2}+\frac{Q(f_1)-Q(f_2)}{2(f_1-f_2)}(g_1+g_2)-g_1g_2.$$
These definitions are equivalent by virtue of (\ref{relod}).
Note that $\frac{H(f_1)-H(f_2)}{f_1-f_2}\in S^2\cal F\subset\cal F\otimes\cal F$ for an arbitrary polynomial $H$. Indeed, $\frac{f_1^m-f_2^m}{f_1-f_2}=f_1^{m-1}+f_1^{m-2}f_2+...+f_2^{m-1}\in S^2\cal F$.

We define elements $x_0,x_2,x_3,x_4,...\in\cal F$ by
\begin{equation}\label{x}
 x_{2i}=f^i,~x_{2i+3}=f^ig,~i=0,1,2,...
\end{equation}
Let ${\cal F}_n\subset \cal F$ be an $n$-dimensional linear space with a basis $\{x_0,x_2,x_3,...,x_n\}=\{x_0,x_i;~2\leq i\leq n\}$. Note that ${\cal F}_n\subset \cal F$ is not a subalgebra 
of $\cal F$. We assume ${\cal F}_n\subset {\cal} F_{ev}$ if $n$ is even and ${\cal F}_n\subset {\cal} F_{od}$ if $n$ is odd. We will identify  $S^*{\cal F}_n$ with a polynomial algebra $\C[x_0,x_2,..,x_n]$ in $n$ variables. 
In particular:
\begin{equation}\label{id}
f_1^if_2^j+f_1^jf_2^i=x_{2i}x_{2j},~f_1^if_2^jg_2+f_1^jf_2^ig_1=x_{2i}x_{2j+3},~(f_1^if_2^j+f_1^jf_2^i)g_1g_2=x_{2i+3}x_{2j+3}. 
\end{equation}

\subsection{Construction in the case of even $n$}

{\bf Proposition 1.} The following formula
\begin{equation}\label{poisev}
\{\phi,\psi\}=n\lambda_{ev}(\phi_1\psi_2-\psi_1\phi_2)+\phi_1D(\psi_2)+\phi_2D(\psi_1)-\psi_1D(\phi_2)-\psi_2D(\phi_1) 
\end{equation}
defines a quadratic Poisson bracket on the polynomial ring $S^*{\cal F}_n=\C[x_0,x_2,..,x_n]$ where $n$ is even. Here $\phi,\psi\in{\cal F}_n$ and $\{\phi,\psi\}\in S^2{\cal F}_n$. This Poisson bracket is linear 
with respect to coefficients $a_0,...,a_4,b_0,...,b_2$ of polynomials $P_{ev},~Q$ and, 
therefore, can be written in the form $\{\cdot,\cdot\}=\{\cdot,\cdot\}_0+\sum_{i=0}^4a_i\{\cdot,\cdot\}_{i,1}+\sum_{j=0}^2b_j\{\cdot,\cdot\}_{j,2}$ where 
$\{\cdot,\cdot\}_0,~\{\cdot,\cdot\}_{i,1},~\{\cdot,\cdot\}_{j,2}$ are pairwise compatible. Therefore, for each even $n$ we have constructed nine compatible quadratic Poisson brackets 
in $n$ variables.

{\bf Proof.} The Jacobi identity is a consequence of a functional
construction described in the next section. Let us check linearity
with respect to coefficients of $P_{ev},~Q$.  Each of $\phi,\psi\in{\cal F}_n\subset{\cal F}_{ev}$
can be of the form $R(f)$ or $R(f)g$ where $R$ is a
polynomial. Therefore, we have three cases:

\noindent
{\bf Case 1.} Let $\phi=R(f),~\psi=T(f)$. We have 
\begin{eqnarray*}
\{\phi,\psi\}&=&\{R(f),T(f)\} \\
&=&n\lambda_{ev}(R(f_1)T(f_2)-R(f_2)T(f_1))+R(f_1)D(T(f_2))+R(f_2)D(T(f_1))\\
& &-T(f_1)D(R(f_2))-T(f_2)D(R(f_1)) \\
&=&n\frac{R(f_1)T(f_2)-R(f_2)T(f_1)}{f_1-f_2}
   \left(g_1+g_2-\frac{1}{2}Q(f_1)-\frac{1}{2}Q(f_2)\right)+ \\
& &(R(f_2)T^{\prime}(f_1)-T(f_2)R^{\prime}(f_1))(2g_1-Q(f_1))+
   (R(f_1)T^{\prime}(f_2)-T(f_1)R^{\prime}(f_2))(2g_2-Q(f_2))
\end{eqnarray*}
{\bf Case 2.} Let $\phi=R(f),~\psi=T(f)g$. We have 
\begin{eqnarray*}
\{\phi,\psi\}&=&\{R(f),T(f)\} \\
&=&n\lambda_{ev}(R(f_1)T(f_2)g_2-R(f_2)T(f_1)g_1)+R(f_1)D(T(f_2)g_2)+R(f_2)D(T(f_1)g_1)\\
& &-T(f_1)g_1D(R(f_2))-T(f_2)g_2D(R(f_1))\\
&=&n\frac{R(f_1)T(f_2)-R(f_2)T(f_1)}{f_1-f_2}g_1g_2
   -\frac{n}{2}\frac{Q(f_1)-Q(f_2)}{f_1-f_2}\left(R(f_1)T(f_2)g_2+R(f_2)T(f_1)g_1\right)+\\
& &n\frac{R(f_1)T(f_2)P_{ev}(f_2)-R(f_2)T(f_1)P_{ev}(f_1)}{f_1-f_2}\\
& &+R(f_1)T(f_2)(P_{ev}^{\prime}(f_2)+Q^{\prime}(f_2)g_2)
   +R(f_2)T(f_1)(P_{ev}^{\prime}(f_1)+Q^{\prime}(f_1)g_1)\\
& &+R^{\prime}(f_1)T(f_2)(Q(f_1)-2g_1)g_2
   +R^{\prime}(f_2)T(f_1)(Q(f_2)-2g_2)g_1\\
& &+R(f_1)T^{\prime}(f_2)(2P_{ev}(f_2)+Q(f_2)g_2)
   +R(f_2)T^{\prime}(f_1)(2P_{ev}(f_1)+Q(f_1)g_1).
\end{eqnarray*}
Here we used $\lambda_{ev}(R(f_1)T(f_2)g_2-R(f_2)T(f_1)g_1)=$
$$(g_1+g_2)\frac{R(f_1)T(f_2)-R(f_2)T(f_1)}{2(f_1-f_2)}\lambda_{ev}(f_1-f_2)
-\frac{1}{2}(R(f_1)T(f_2)+R(f_2)T(f_1))\lambda_{ev}(g_1-g_2).$$
{\bf Case 3.} Let $\phi=R(f)g,~\psi=T(f)g$. We have 
\begin{eqnarray*}
\{\phi,\psi\}&=&\{R(f),T(f)\} \\
&=&n\lambda_{ev}((R(f_1)T(f_2)-R(f_2)T(f_1))g_1g_2)+R(f_1)g_1D(T(f_2)g_2)+R(f_2)g_2D(T(f_1)g_1)\\
& &-T(f_1)g_1D(R(f_2)g_2)-T(f_2)g_2D(R(f_1)g_1) \\
&=&n\frac{R(f_1)T(f_2)-R(f_2)T(f_1)}{f_1-f_2}(g_1^2g_2+g_1g_2^2
   -\frac{1}{2}Q(f_1)g_1g_2-\frac{1}{2}Q(f_2)g_1g_2)\\
& &+(R(f_2)T^{\prime}(f_1)-T(f_2)R^{\prime}(f_1))(2g_1-Q(f_1))g_1g_2\\
& &+(R(f_1)T^{\prime}(f_2)-T(f_1)R^{\prime}(f_2))(2g_2-Q(f_2))g_1g_2\\
& &+(R(f_1)T(f_2)-R(f_2)T(f_1))(g_1D(g_2)-D(g_1)g_2) \\
&=&n\frac{R(f_1)T(f_2)-R(f_2)T(f_1)}{f_1-f_2}(P_{ev}(f_1)g_2+g_1P_{ev}(f_2)
   +\frac{1}{2}Q(f_1)g_1g_2+\frac{1}{2}Q(f_2)g_1g_2)\\
& &+(R(f_2)T^{\prime}(f_1)-T(f_2)R^{\prime}(f_1))(2P_{ev}(f_1)+Q(f_1)g_1)g_2\\
& &+(R(f_1)T^{\prime}(f_2)-T(f_1)R^{\prime}(f_2))(2P_{ev}(f_2)+Q(f_2)g_2)g_1\\
& &+(R(f_1)T(f_2)-R(f_2)T(f_1))(P_{ev}^{\prime}(f_2)g_1+Q^{\prime}(f_2)g_1g_2
   -P_{ev}^{\prime}(f_1)g_2-Q^{\prime}(f_1)g_1g_2)
\end{eqnarray*}
In these computations we use formulas (\ref{relev}) and (\ref{derev}) 
where $f,g$ are replaced by $f_1,g_1$ or $f_2,g_2$.
Note that in each case we obtain an expression for $\{\phi,\psi\}$
that is linear non-homogeneous in $P_{ev},~Q$ and bi-linear non-homogeneous in
$g_1,~g_2$. Using identifications (\ref{id}) we can write each of
these expressions as a quadratic homogeneous polynomial in
$x_0,x_2,x_3,...,x_n$ with coefficients linear in
$a_0,...,a_4,b_0,...,b_2$.

{\bf Remark 1.} Let us describe an algorithm for the computation of
$\{x_i,x_j\}$ as a quadratic polynomial in $x_0,x_2,x_3,...,x_n$.  One
uses the formulas for $\{\phi,\psi\}$ from the proof of the
proposition 1 where $R(f)=f^i,~T(f)=f^j$ and $P_{ev},~Q$ are given by
(\ref{pol}). For the computation of $\{x_{2i},x_{2j}\}$ the formula
of case 1 is used, for the computation of $\{x_{2i},x_{2j+3}\}$ the
formula of case 2 is used, and for the computation of
$\{x_{2i+3},x_{2j+3}\}$ the formula of case 3 is used.  These formulas
give a polynomial in $f_1,~f_2,~g_1,~g_2$ linear in $g_1,~g_2$. This
polynomial is symmetric under transformations $f_1\leftrightarrow
f_2,~g_1\leftrightarrow g_2$. By using identifications (\ref{id}) the
expressions can be written as a polynomial quadratic in
$x_0,x_2,x_3,...,x_n$.

{\bf Remark 2.} In all cases of the above remark $\{x_i,x_j\}$ is
linear non-homogeneous with respect to the eight coefficients $a_0,...,b_2$ of
polynomials $P_{ev},~Q$. Therefore, nine compatible Poisson brackets can be
obtained in the following way:
$$\{x_i,x_j\}_0=\{x_i,x_j\}|_{a_0=...=b_2=0},
 ~\{x_i,x_j\}_{k,1}=\frac{\partial\{x_i,x_j\}}{\partial a_k},
 ~\{x_i,x_j\}_{k,2}=\frac{\partial\{x_i,x_j\}}{\partial b_k}.$$

\subsection{Construction in the case of odd $n$}

{\bf Proposition 1$^{\prime}$.} The following formula
\begin{equation}\label{poisod}
\{\phi,\psi\}=(n\lambda_{od}+g_2-g_1+\frac{1}{2}b_2(f_1-f_2))(\phi_1\psi_2-\psi_1\phi_2)+\phi_1D(\psi_2)+\phi_2D(\psi_1)-\psi_1D(\phi_2)-\psi_2D(\phi_1) 
\end{equation}
defines a quadratic Poisson bracket on the polynomial ring $S^*{\cal F}_n=\C[x_0,x_2,..,x_n]$ where $n$ is odd. Here $\phi,\psi\in{\cal F}_n$ and $\{\phi,\psi\}\in S^2{\cal F}_n$. This Poisson bracket is linear 
with respect to $c$ and coefficients $a_0,...,a_3,b_0,...,b_2$ of polynomials $P_{od},~Q$ and, 
therefore, can be written in the form $\{\cdot,\cdot\}=\{\cdot,\cdot\}_0+\sum_{i=0}^3a_i\{\cdot,\cdot\}_{i,1}+\sum_{j=0}^2b_j\{\cdot,\cdot\}_{j,2}+c\{\cdot,\cdot\}_{3}$ where 
$\{\cdot,\cdot\}_0,~\{\cdot,\cdot\}_{i,1},~\{\cdot,\cdot\}_{j,2},~\{\cdot,\cdot\}_{3}$ are pairwise compatible. Therefore, for each odd $n$ we have constructed nine compatible quadratic Poisson brackets 
in $n$ variables.

{\bf Proof.} The Jacobi identity is a consequence of a functional
construction described in the next section. Let us check linearity
with respect to $c$ and coefficients of $P_{od},~Q$.  Each of $\phi,\psi\in{\cal F}_n\subset{\cal F}_{ev}$
can be of the form $R(f)$ or $R(f)g$ where $R$ is a
polynomial. Therefore, we have three cases:

\noindent
{\bf Case 1.} Let $\phi=R(f),~\psi=T(f)$. We have 
\begin{eqnarray*}
\{\phi,\psi\}&=&\{R(f),T(f)\} \\
&=&(n\lambda_{od}+g_2-g_1+\frac{1}{2}b_2(f_1-f_2))(R(f_1)T(f_2)-R(f_2)T(f_1))+\\& &R(f_1)D(T(f_2))+R(f_2)D(T(f_1))\\
& &-T(f_1)D(R(f_2))-T(f_2)D(R(f_1)) \\
&=&n\frac{R(f_1)T(f_2)-R(f_2)T(f_1)}{f_1-f_2}
   \left((f_1+c)g_1+(f_2+c)g_2-\frac{1}{2}Q(f_1)-\frac{1}{2}Q(f_2)\right)+\\& &(g_2-g_1+\frac{1}{2}b_2(f_1-f_2))(R(f_1)T(f_2)-R(f_2)T(f_1))+ \\
& &(R(f_2)T^{\prime}(f_1)-T(f_2)R^{\prime}(f_1))(2(f_1+c)g_1-Q(f_1))+\\& &
   (R(f_1)T^{\prime}(f_2)-T(f_1)R^{\prime}(f_2))(2(f_2+c)g_2-Q(f_2))
\end{eqnarray*}
{\bf Case 2.} Let $\phi=R(f),~\psi=T(f)g$. We have 
\begin{eqnarray*}
\{\phi,\psi\}&=&\{R(f),T(f)\} \\
&=&(n\lambda_{od}+g_2-g_1+\frac{1}{2}b_2(f_1-f_2))(R(f_1)T(f_2)g_2-R(f_2)T(f_1)g_1)+\\& &R(f_1)D(T(f_2)g_2)+R(f_2)D(T(f_1)g_1)\\
& &-T(f_1)g_1D(R(f_2))-T(f_2)g_2D(R(f_1))\\
&=&\frac{n}{2}\frac{R(f_1)T(f_2)-R(f_2)T(f_1)}{f_1-f_2}(f_1+f_2+2c)g_1g_2
   -\\& &\frac{n}{2}\frac{Q(f_1)-Q(f_2)}{f_1-f_2}\left(R(f_1)T(f_2)g_2+R(f_2)T(f_1)g_1\right)+\\& &\frac{n-2}{2}(R(f_1)T(f_2)+R(f_2)T(f_1))g_1g_2+\\
& &n\frac{R(f_1)T(f_2)P_{od}(f_2)-R(f_2)T(f_1)P_{od}(f_1)}{f_1-f_2}+\\& &\frac{1}{2}b_2(f_1-f_2)(R(f_1)T(f_2)g_2-R(f_2)T(f_1)g_1)\\
& &+R(f_1)T(f_2)(P_{od}^{\prime}(f_2)+Q^{\prime}(f_2)g_2)
   +R(f_2)T(f_1)(P_{od}^{\prime}(f_1)+Q^{\prime}(f_1)g_1)\\
& &+R^{\prime}(f_1)T(f_2)(Q(f_1)-2(f_1+c)g_1)g_2
   +R^{\prime}(f_2)T(f_1)(Q(f_2)-2(f_2+c)g_2)g_1\\
& &+R(f_1)T^{\prime}(f_2)(2P_{od}(f_2)+Q(f_2)g_2)
   +R(f_2)T^{\prime}(f_1)(2P_{od}(f_1)+Q(f_1)g_1).
\end{eqnarray*}
Here we used $\lambda_{od}(R(f_1)T(f_2)g_2-R(f_2)T(f_1)g_1)=$
$$(g_1+g_2)\frac{R(f_1)T(f_2)-R(f_2)T(f_1)}{2(f_1-f_2)}\lambda_{od}(f_1-f_2)
-\frac{1}{2}(R(f_1)T(f_2)+R(f_2)T(f_1))\lambda_{od}(g_1-g_2).$$
{\bf Case 3.} Let $\phi=R(f)g,~\psi=T(f)g$. We have 
\begin{eqnarray*}
\{\phi,\psi\}&=&\{R(f),T(f)\} \\
&=&(n\lambda_{od}+g_2-g_1+\frac{1}{2}b_2(f_1-f_2))((R(f_1)T(f_2)-R(f_2)T(f_1))g_1g_2)+\\& &R(f_1)g_1D(T(f_2)g_2)+R(f_2)g_2D(T(f_1)g_1)\\
& &-T(f_1)g_1D(R(f_2)g_2)-T(f_2)g_2D(R(f_1)g_1) \\
&=&n\frac{R(f_1)T(f_2)-R(f_2)T(f_1)}{f_1-f_2}((f_1+c)g_1^2g_2+(f_2+c)g_1g_2^2
   -\frac{1}{2}Q(f_1)g_1g_2-\frac{1}{2}Q(f_2)g_1g_2)+\\& &(R(f_1)T(f_2)-R(f_2)T(f_1))(\frac{1}{2}b_2(f_1-f_2)g_1g_1+g_1g_2^2-g_1^2g_2)\\
& &+(R(f_2)T^{\prime}(f_1)-T(f_2)R^{\prime}(f_1))(2(f_1+c)g_1-Q(f_1))g_1g_2\\
& &+(R(f_1)T^{\prime}(f_2)-T(f_1)R^{\prime}(f_2))(2(f_2+c)g_2-Q(f_2))g_1g_2\\
& &+(R(f_1)T(f_2)-R(f_2)T(f_1))(g_1D(g_2)-D(g_1)g_2) \\
&=&n\frac{R(f_1)T(f_2)-R(f_2)T(f_1)}{f_1-f_2}(P_{od}(f_1)g_2+g_1P_{od}(f_2)
   +\frac{1}{2}Q(f_1)g_1g_2+\frac{1}{2}Q(f_2)g_1g_2)\\& &+\frac{b_2}{2}(R(f_1)T(f_2)-R(f_2)T(f_1))(f_1-f_2)g_1g_2\\
& &+(R(f_2)T^{\prime}(f_1)-T(f_2)R^{\prime}(f_1))(2P_{od}(f_1)+Q(f_1)g_1)g_2\\
& &+(R(f_1)T^{\prime}(f_2)-T(f_1)R^{\prime}(f_2))(2P_{od}(f_2)+Q(f_2)g_2)g_1\\
& &+(R(f_1)T(f_2)-R(f_2)T(f_1))(P_{od}^{\prime}(f_2)g_1+Q^{\prime}(f_2)g_1g_2
   -P_{od}^{\prime}(f_1)g_2-Q^{\prime}(f_1)g_1g_2)
\end{eqnarray*}
In these computations we use formulas (\ref{relod}) and (\ref{derod}) where $f,g$ are replaced by $f_1,g_1$ or $f_2,g_2$.
Note that in each case we obtain an expression for $\{\phi,\psi\}$
that is linear non-homogeneous in $P_{od},~Q$ and bi-linear non-homogeneous in
$g_1,~g_2$. Using identifications (\ref{id}) we can write each of
these expressions as a quadratic homogeneous polynomial in
$x_0,x_2,x_3,...,x_n$ with coefficients linear in
$a_0,...,a_3,b_0,...,b_2,c$.

{\bf Remark 1$^\prime$.} Let us describe an algorithm for the computation of
$\{x_i,x_j\}$ as a quadratic polynomial in $x_0,x_2,x_3,...,x_n$.  One
uses the formulas for $\{\phi,\psi\}$ from the proof of the
proposition 1$^\prime$ where $R(f)=f^i,~T(f)=f^j$ and $P_{od},~Q$ are given by
(\ref{pol}). For the computation of $\{x_{2i},x_{2j}\}$ the formula
of case 1 is used, for the computation of $\{x_{2i},x_{2j+3}\}$ the
formula of case 2 is used, and for the computation of
$\{x_{2i+3},x_{2j+3}\}$ the formula of case 3 is used.  These formulas
give a polynomial in $f_1,~f_2,~g_1,~g_2$ linear in $g_1,~g_2$. This
polynomial is symmetric under transformations $f_1\leftrightarrow
f_2,~g_1\leftrightarrow g_2$. By using identifications (\ref{id}) the
expressions can be written as a polynomial quadratic in
$x_0,x_2,x_3,...,x_n$.

{\bf Remark 2$^\prime$.} In all cases of the above remark $\{x_i,x_j\}$ is
linear non-homogeneous with respect to $c$ and the seven coefficients $a_0,...,b_2$ of
polynomials $P_{od},~Q$. Therefore, nine compatible Poisson brackets can be
obtained in the following way:
$$\{x_i,x_j\}_0=\{x_i,x_j\}|_{a_0=...=b_2=c=0},
 ~\{x_i,x_j\}_{k,1}=\frac{\partial\{x_i,x_j\}}{\partial a_k},
 ~\{x_i,x_j\}_{k,2}=\frac{\partial\{x_i,x_j\}}{\partial b_k},
 ~\{x_i,x_j\}_{3}=\frac{\partial\{x_i,x_j\}}{\partial c}.$$

{\bf Remark 3.} If we replace $n$ in the formulas (\ref{poisev}), (\ref{poisod}) by an arbitrary constant $\alpha$, then these formulas 
still define Poisson brackets on the polynomial algebra $\C[x_0,x_2,x_3,...]$. However, $\C[x_0,x_2,x_3,...,x_n]$ $\subset$ $\C[x_0,x_2,x_3,...]$ 
is closed with respect to these brackets only if $\alpha=n$.

\section{Functional construction} \label{functional}

\subsection{General constructions}

Recall a general construction of associative algebras and Poisson
structures \cite{odes}. Let $\lambda(x,y)$ be a meromorphic function
in two variables. We construct an associative algebra $A_{\lambda}$
through:
$$A_{\lambda}=\C\oplus F_1\oplus F_2\oplus F_3\oplus...$$ where $F_m$
is the space of symmetric meromorphic functions in $m$ variables and a
product $f\star g\in F_{a+b}$ of $f\in F_a,~g\in F_b$ is defined by:
$$f\star g(z_1,...,z_{a+b})=\frac{1}{a!b!}\sum_{\sigma\in
  S_{a+b}}f(z_{\sigma_1},...,z_{\sigma_a})g(z_{\sigma_{a+1}},...,z_{\sigma_{a+b}})\prod_{1\leq
  p\leq a,~ a+1\leq q\leq a+b}\lambda(z_{\sigma_p},z_{\sigma_q}).$$
Note that this formula defines an associative product for an arbitrary
function $\lambda$ and this product is non-commutative if $\lambda$ is
not symmetric. In particular, if $\lambda(x,y)=1+
\frac{1}{2}\epsilon\mu(x,y)+o(\epsilon)$, then we obtain a Poisson
algebra. Assume $\mu(x,y)=-\mu(y,x)$. For $f,~g\in F_1$ we get the
following formulas for the associative commutative product $fg\in F_2$
and the Poisson bracket:
$$fg(x,y)=f(x)g(y)+g(x)f(y),~\{f,g\}=(f(x)g(y)-g(x)f(y))\mu(x,y).$$ If
$\lambda=\lambda(x-y)$, then the formula for the product can be
deformed in the following way:
$$f\star g(z_1,...,z_{a+b})=\frac{1}{a!b!}\sum_{\sigma\in
  S_{a+b}}f(z_{\sigma_1},...,z_{\sigma_a})g(z_{\sigma_{a+1}}+ap,...,z_{\sigma_{a+b}}+ap)\prod_{1\leq
  p\leq a,~ a+1\leq q\leq a+b}\lambda(z_{\sigma_p}-z_{\sigma_q})$$
where $p$ is an arbitrary constant. The corresponding formula for the
Poisson brackets (if we set $p=\epsilon\alpha$) is:
$$\{f,g\}=(f(x)g(y)-g(x)f(y))\mu(x-y)
          +\alpha(f(x)g^{\prime}(y)+f(y)g^{\prime}(x)-g(x)f^{\prime}(y)-g(y)f^{\prime}(x)).$$
We will need the following generalization of the last formula:
\begin{eqnarray}
\{f,g\}&=&(f(x)g(y)-g(x)f(y))(\mu(x-y)+\nu(x)-\nu(y)) \nonumber \\ 
       & &+\alpha(f(x)g^{\prime}(y)+f(y)g^{\prime}(x)
          -g(x)f^{\prime}(y)-g(y)f^{\prime}(x)).  \label{poisgen}
\end{eqnarray}
Here $\nu$ is an arbitrary function. This formula is obtained from the
previous one by transformation $f\to\kappa f,~g\to\kappa
g,~\{f,g\}\to\kappa^2\{f,g\}$ where $\alpha\kappa^{\prime}=-\nu$.

Note that the algebra $A_{\lambda}$ is very large. One can construct
associative algebras (and the corresponding Poisson algebras) of a
reasonable size by a suitable choice of spaces $F_{\alpha}$ and
function $\lambda$. For example, the algebra $Q_n(\tau,\eta)$ and the
Poisson algebra $q_n(\tau)$ can be constructed in this way
\cite{odes}. See \cite{odfeig} for the functional construction of 
a wider class of Poisson algebras.

\subsection{The case of even $n$}

It is clear that equation (\ref{relev}) defines an elliptic curve in
$\C^2$ with coordinates $f,~g$. Therefore, one can find elliptic
functions $f=f(z),~g=g(z)$ such that
\begin{equation}\label{rel1ev}
g(z)^2=P_{ev}(f(z))+Q(f(z))g(z). 
\end{equation}
Moreover, one can assume (see (\ref{derev}))
\begin{equation}\label{der1ev}
   f^{\prime}(z)=2g(z)-Q(f(z)),
~~~g^{\prime}(z)=P_{ev}^{\prime}(f(z))+Q^{\prime}(f(z))g(z).
\end{equation}
Note that elliptic functions $f(z),~g(z)$ have a form:
$$  f(z)=c_1+c_2\zeta(z,\tau)+c_3\zeta(z-u,\tau),
 ~~~g(z)=c_4+c_5\zeta(z,\tau)+c_6\zeta(z-u,\tau)
         +c_7\zeta^{\prime}(z,\tau)+c_8\zeta^{\prime}(z-u,\tau).$$
Here $\zeta(z,\tau)$ is the Weierstrass elliptic function, $\tau$ is a
modular parameter and constants $c_1,...,c_8, u, \tau$ are determined
by $a_0,...,b_2$.  There exists an elliptic function in two variables
$\mu_{ev}(z_1,z_2)$ such that
\begin{equation}\label{hom1ev}
(f(z_1)-f(z_2))\mu_{ev}(z_1,z_2)=g(z_1)+g(z_2)-\frac{1}{2}Q(f(z_1))-\frac{1}{2}Q(f(z_2)),
\end{equation}
$$(g(z_1)-g(z_2))\mu_{ev}(z_1,z_2)=\frac{P_{ev}(f(z_1))-P_{ev}(f(z_2))}{f(z_1)-f(z_2)}+\frac{Q(f(z_1))-Q(f(z_2))}{2(f(z_1)-f(z_2))}(g(z_1)+g(z_2)).$$
This function has the form 
$$\mu_{ev}(z_1,z_2)=\zeta(z_1-z_2,\tau)+k_1\zeta(z_1,\tau)+k_2\zeta(z_1-u,\tau)-k_1\zeta(z_2,\tau)-k_2\zeta(z_2-u,\tau)$$ 
for some constants $k_1,~k_2$ such that $k_1+k_2=1$.

Let ${\cal F}_n$ be the space of elliptic functions in one variable with
periods 1 and $\tau$, holomorphic outside $z=0,~u$ modulo periods and having poles of order not larger than $n$ at $z=0,~u$.
It is clear that $\{e_0(z),e_2(z),e_3(z),e_4(z),...,e_n(z)\}$ is a basis of
the linear space ${\cal F}_n$ where we define
$$e_{2i}(z)=f(z)^i,~e_{2i+3}(z)=f(z)^ig(z),~i=0,1,2,...  .$$ We will
identify $S^m{\cal F}_n$ with the space of symmetric elliptic functions in
$m$ variables $\{h(z_1,...,z_m)\}$ holomorphic if $z_k\ne 0,u$ modulo
periods and having poles of order not larger than $n$ at $z_k=0,~u$.  We construct a bilinear operator $\{,\} : \Lambda^2{\cal F}_n\to
S^2{\cal F}_n$ as follows: for $\phi,\psi\in{\cal F}_n$ we set
\begin{eqnarray}
\{\phi,\psi\}(z_1,z_2)&=&n\mu_{ev}(z_1,z_2)(\phi(z_1)\psi(z_2)-\psi(z_1)\phi(z_2)) 
 \nonumber \\
& &+\phi(z_1)\psi^{\prime}(z_2)
   +\phi(z_2)\psi^{\prime}(z_1)
   -\psi(z_1)\phi^{\prime}(z_2)
   -\psi(z_2)\phi^{\prime}(z_1).  \label{pois1ev}
\end{eqnarray}

{\bf Proposition 2.} The formula (\ref{pois1ev}) defines a Poisson
structure on the polynomial algebra $S^{*}{\cal F}_n$. This Poisson
bracket is linear with respect to coefficients
$a_0,...,a_4,b_0,...,b_2$ of polynomials $P_{ev},~Q$ and, therefore, can be
written in the form
$\{\cdot,\cdot\}=                \{\cdot,\cdot\}_0
                  +\sum_{i=0}^4a_i\{\cdot,\cdot\}_{i,1}
                  +\sum_{j=0}^2b_j\{\cdot,\cdot\}_{j,2}$
where
$\{\cdot,\cdot\}_0,~\{\cdot,\cdot\}_{i,1},~\{\cdot,\cdot\}_{j,2}$ are
pairwise compatible. Therefore, for each even $n$ we have constructed nine compatible
quadratic Poisson brackets in $n$ variables.

{\bf Proof.} This is just a reformulation of the Proposition 1. Formula (\ref{pois1ev}) is a special case of (\ref{poisgen}) 
and therefore the Jacobi identity for (\ref{pois1ev}) is satisfied. One can check straightforwardly that if $\phi,\psi\in{\cal F}_n$, then $\{\phi,\psi\}(z_1,z_2)$ given by 
(\ref{pois1ev}) is a symmetric elliptic function in two variables $z_1,z_2$ having poles of order not larger than $n$ at $z_1,z_2=0,u$ and therefore $\{\phi,\psi\}(z_1,z_2)\in S^2{\cal F}_n$.

\subsection{The case of odd $n$}

It is clear that equation (\ref{relod}) defines an elliptic curve in
$\C^2$ with coordinates $f,~g$. Therefore, one can find elliptic
functions $f=f(z),~g=g(z)$ such that
\begin{equation}\label{rel1od}
(f(z)+c)g(z)^2=P_{od}(f(z))+Q(f(z))g(z). 
\end{equation}
Moreover, one can assume (see (\ref{derod}))
\begin{equation}\label{der1od}
   f^{\prime}(z)=2(f(z)+c)g(z)-Q(f(z)),
~~~g^{\prime}(z)=P_{od}^{\prime}(f(z))+Q^{\prime}(f(z))g(z)-g(z)^2.
\end{equation}
Note that elliptic functions $f(z),~g(z)$ have a form:
$$  f(z)=c_1+c_2\zeta(z,\tau)+c_3\zeta(z-u,\tau),
 ~~~g(z)=c_4+c_5\zeta(z,\tau)+c_6\zeta(z-u,\tau)+c_7\zeta(z-v,\tau).$$
Here $\zeta(z,\tau)$ is the Weierstrass elliptic function, $\tau$ is a
modular parameter and constants $c_1,...,c_7, u, v, \tau$ are determined
by $a_0,...,b_2,c$.  There exists an elliptic function in two variables
$\mu_{od}(z_1,z_2)$ such that
\begin{equation}\label{hom1od}
(f(z_1)-f(z_2))\mu_{od}(z_1,z_2)=(f(z_1)+c)g(z_1)+(f(z_2)+c)g(z_2)-\frac{1}{2}Q(f(z_1))-\frac{1}{2}Q(f(z_2)),
\end{equation}
$$(g(z_1)-g(z_2))\mu_{od}(z_1,z_2)=\frac{P(f(z_1))-P(f(z_2))}{f(z_1)-f(z_2)}+\frac{Q(f(z_1))-Q(f(z_2))}{2(f(z_1)-f(z_2))}(g(z_1)+g(z_2))-g(z_1)g(z_2).$$
This function has the form 
$$\mu_{od}(z_1,z_2)=\zeta(z_1-z_2,\tau)+k_1\zeta(z_1,\tau)+k_2\zeta(z_1-u,\tau)-k_1\zeta(z_2,\tau)-k_2\zeta(z_2-u,\tau)$$ 
for some constants $k_1,~k_2$ such that $k_1+k_2=1$.

Let ${\cal F}_n$ be the space of elliptic functions in one variable with
periods 1 and $\tau$, holomorphic outside $z=0,~u,~v$ modulo periods and having poles of order not larger than $n$ at $z=0,~u$ and not larger than one at $z=v$.
It is clear that $\{e_0(z),e_2(z),e_3(z),e_4(z),...,e_n(z)\}$ is a basis of
the linear space ${\cal F}_n$ where we define
$$e_{2i}(z)=f(z)^i,~e_{2i+3}(z)=f(z)^ig(z),~i=0,1,2,...  .$$ We will
identify $S^m{\cal F}_n$ with the space of symmetric elliptic functions in
$m$ variables $\{h(z_1,...,z_m)\}$ holomorphic if $z_k\ne 0,u,v$ modulo
periods and having poles of order not larger than $n$ at $z_k=0,~u$ and not larger than one at $z_k=v$.  We construct a bilinear operator $\{,\} : \Lambda^2{\cal F}_n\to
S^2{\cal F}_n$ as follows: for $\phi,\psi\in{\cal F}_n$ we set
\begin{eqnarray}
\{\phi,\psi\}(z_1,z_2)&=&(n\mu_{od}(z_1,z_2)+g(z_2)-g(z_1)+\frac{1}{2}b_2(f(z_1)-f(z_2)))(\phi(z_1)\psi(z_2)-\psi(z_1)\phi(z_2)) 
 \nonumber \\
& &+\phi(z_1)\psi^{\prime}(z_2)
   +\phi(z_2)\psi^{\prime}(z_1)
   -\psi(z_1)\phi^{\prime}(z_2)
   -\psi(z_2)\phi^{\prime}(z_1).  \label{pois1od}
\end{eqnarray}

{\bf Proposition 2$^{\prime}$.} The formula (\ref{pois1od}) defines a Poisson
structure on the polynomial algebra $S^{*}{\cal F}_n$. This Poisson
bracket is linear with respect to $c$ and coefficients
$a_0,...,a_3,b_0,...,b_2$ of polynomials $P_{od},~Q$ and, therefore, can be
written in the form
$\{\cdot,\cdot\}=                \{\cdot,\cdot\}_0
                  +\sum_{i=0}^3a_i\{\cdot,\cdot\}_{i,1}
                  +\sum_{j=0}^2b_j\{\cdot,\cdot\}_{j,2}+c\{\cdot,\cdot\}_3$
where
$\{\cdot,\cdot\}_0,~\{\cdot,\cdot\}_{i,1},~\{\cdot,\cdot\}_{j,2},~\{\cdot,\cdot\}_3$ are
pairwise compatible. Therefore, for each odd $n$ we have constructed nine compatible
quadratic Poisson brackets in $n$ variables.

{\bf Proof.} This is just a reformulation of the Proposition 1$^{\prime}$. Formula (\ref{pois1od}) is a special case of (\ref{poisgen}) 
and therefore the Jacobi identity for (\ref{pois1od}) is satisfied. One can check straightforwardly that if $\phi,\psi\in{\cal F}_n$, then $\{\phi,\psi\}(z_1,z_2)$ given by 
(\ref{pois1od}) is a symmetric elliptic function in two variables $z_1,z_2$ having poles of order not larger than $n$ at $z_1,z_2=0,u$ and not larger than one at $z_1,z_2=v$ and 
therefore $\{\phi,\psi\}(z_1,z_2)\in S^2{\cal F}_n$.

\section{Symplectic leaves and Casimir elements} \label{cas}

For $p,n\in\N$ we denote by $b_{p,n}$ the Poisson algebra spanned by the elements
$$\{h(f_1,g_1,...,f_p,g_p)e_1^{\alpha_1}...e_p^{\alpha_p};\alpha_1,...,\alpha_p\in\Z_{\geq 0}\}$$ 
as a linear space, where $h$ is a rational function,
$f_i,g_i$ for each $i=1,...,p$ are subject to relation (\ref{relev})
if $n$ is even and (\ref{relod}) if $n$ is odd where $f,g$ are
replaced by $f_i,g_i$. In other words\footnote{Recall that ${\cal
    F}={\cal F}_{ev}$ if $n$ is even and ${\cal F}={\cal F}_{od}$ if
  $n$ is odd.}, $h\in Frac(\otimes^p\cal F)$. A Poisson bracket on
$b_{p,n}$ is defined as follows:
$$\{f_i,e_j\}=-D(f_i)e_j,~\{g_i,e_j\}=-D(g_i)e_j,~\{f_i,e_i\}=\frac{n-2}{2}D(f_i)e_j,~\{g_i,e_i\}=\frac{n-2}{2}D(g_i)e_i$$
where $i\ne j$, $D$ is defined by (\ref{derev}) for even $n$ and
(\ref{derod}) for odd $n$. We also assume $$\{e_i,e_j\}=n\lambda_{i,j,ev}e_ie_j$$
if $n$ is even and 
$$\{e_i,e_j\}=\left(n\lambda_{i,j,od}+\frac{1}{2}b_2(f_i-f_j)+g_j-g_i\right)e_ie_j$$
if $n$ is odd. Here $\lambda_{i,j,ev}$ (resp. $\lambda_{i,j,od}$) is given by (\ref{homev}) 
(resp. (\ref{homod})) where $f_1,g_1,f_2,g_2$ are replaced by $f_i,g_i,f_j,g_j$ 
correspondingly. All brackets between $f_i,g_j$ are zero.

Let us define a linear map $\phi_p:{\cal F}_n\to b_{p,n}$ by the formula
$$\phi_p(x_{2j})=\sum_{i=1}^pf^j_ie_i,~~~\phi_p(x_{2j+3})=\sum_{i=1}^pf^j_ig_ie_i.$$
There is a unique extension of this map to the homomorphism
of commutative algebras $S^{*}({\cal F}_n)\to b_{p,n}$ which we also
denote by $\phi_p$.

{\bf Proposition 3.} The map $\phi_p:S^{*}({\cal F}_n)\to b_{p,n}$ is a homomorphism of Poisson algebras.

{\bf Proof.} One can check straightforwardly that
$\phi_p(\{r,s\})=\sum_{i,j=1}^p\{r,s\}_{i,j}e_ie_j$ where $r,s$ are
arbitrary elements from ${\cal F}_n$ and $\{r,s\}_{i,j}$ is obtained
from $\{r,s\}\in {\cal F}_n\otimes {\cal F}_n$ replacing
$f_1,g_1,f_2,g_2$ by $f_i,g_i,f_j,g_j$ correspondingly. This implies
the proposition.

It is known \cite{odes, od1} that if $2p<n$, then the map $\phi_p$
defines a $2p$ dimensional symplectic leaf of the Poisson algebra
$S^{*}({\cal F}_n)$. Moreover, central elements of the Poisson algebra
$S^{*}({\cal F}_n)$ belong to $\ker\phi_p$ for $2p<n$. One can check
that for $p=\frac{n}{2}-1$ for even $n$ (resp. $p=\frac{n-1}{2}$ for
odd $n$) the ideal $\ker\phi_p$ is generated by two elements of degree
$\frac{n}{2}$ (resp. by one element of degree $n$). We denote these
elements by $C_0^{\frac{n}{2}},C_1^{\frac{n}{2}}$ if $n$ is even
(resp. by $C^n$ if $n$ is odd). The center of the Poisson algebra
$S^{*}({\cal F}_n)$ is a polynomial algebra generated by
$C_0^{\frac{n}{2}},C_1^{\frac{n}{2}}$ if $n$ is even (resp. by $C^{n}$
if $n$ is odd).  Let us describe these elements explicitly (see also
\cite{od1}).

{\em Let $n$ be even.} We define elements $x_i,~i\in\Z$ by the formula (\ref{x}). It is 
clear that $x_i\in Frac({\cal F}_{ev})$ if $i=1,-1,-2,...$ and  $x_0,x_2,...,x_n\in{\cal F}_n$. 
If $\frac{n}{2}$ is even we set\footnote{These formulas work for $n>4$. If $n=4$ we set \begin{eqnarray*} 
C_0^{2}&\!\!=\!\!&\det\left((x_0,x_2)^t(x_0,x_2)\right),\\ 
C_1^{2}&\!\!=\!\!&\det\left((x_0,x_1)^t(x_2,x_3)\right)-\det\left((x_0,x_2)^t(x_2,a_4x_{4}+b_2x_{3})\right)-\det\left((a_0x_{-2}+b_0x_1,x_0)^t(x_0,x_2)\right)
\end{eqnarray*}} 
\begin{eqnarray*} 
C_0^{\frac{n}{2}}&\!\!=\!\!&\det\left((x_0,x_2,..,x_{\frac{n}{2}})^t(x_0,x_2,..,x_{\frac{n}{2}})\right),\\ 
C_1^{\frac{n}{2}}&\!\!=\!\!&\ \ \ \det\left((x_0,x_1,x_2,..,x_{\frac{n}{2}-1})^t(x_2,..,x_{\frac{n}{2}+1})\right)\\
&&-\det\left((x_0,x_1,x_2,..,x_{\frac{n}{2}-2},x_{\frac{n}{2}})^t(x_2,..,x_{\frac{n}{2}},a_4x_{\frac{n}{2}+2}+b_2x_{\frac{n}{2}+1})\right)\\
&&-\det\left((a_0x_{-2}+b_0x_1,x_0,x_2,..,x_{\frac{n}{2}-1})^t(x_0,x_2,x_4,..,x_{\frac{n}{2}+1})\right) \\
&&+\det \left( (a_0x_{-2}+b_0x_1,x_0,x_2,..,x_{\frac{n}{2}-2},x_{\frac{n}{2}})^t (x_0,x_2,x_4,..,x_{\frac{n}{2}},a_4x_{\frac{n}{2}+2}+b_2x_{\frac{n}{2}+1}) \right)
\end{eqnarray*}
and if $\frac{n}{2}$ is odd
\begin{eqnarray*}
C_0^{\frac{n}{2}}&\!\!\!=\!\!\!&\det\left((x_0,x_2,..,x_{\frac{n}{2}})^t(x_0,x_2,..,x_{\frac{n}{2}})\right)
-\det\left((x_0,x_2,..,x_{\frac{n}{2}-1},x_{\frac{n}{2}+1})^t(x_0,x_2,..,x_{\frac{n}{2}-1},b_2x_{\frac{n}{2}}+a_4x_{\frac{n}{2}+1})\right),\\ 
C_1^{\frac{n}{2}}&\!\!\!=\!\!\!&\det\left((x_2,..,x_{\frac{n}{2}+1})^t(x_0,x_1,x_2,..,x_{\frac{n}{2}-1})\right)
-\det\left((x_{-2},x_0,x_2,..,x_{\frac{n}{2}-1})^t(a_0x_0+b_0x_3,x_2,x_4,..,x_{\frac{n}{2}+1})\right) .
\end{eqnarray*}
In these formulas we use the product in the algebra $Frac({\cal F}_{ev})$ 
for 
computing 
products of vector components. 
Therefore entries of our
matrices of the form $v^tw$ are linear combinations of $x_i\in
Frac({\cal F}_{ev})$, $i\in\Z$. On the other hand, 
for 
computing
determinants we use the product in $S^{*}(Frac({\cal F}_{ev}))$. So
these determinants are polynomials of degree $\frac{n}{2}$ in
$x_i,~i\in\Z$.\footnote{For example, $\det\left((x_0,x_2)^t(x_0,x_2)\right)=\det\left((1,f)^t(1,f)\right)=\det\left(
\begin{array}{cc}
1 & f  \\
f & f^2   \end{array}
\right)=\det\left(
\begin{array}{cc}
x_0 & x_2  \\
x_2 & x_4   \end{array}
\right)=x_0x_4-x_2^2.$} Moreover, it turns out that in our linear combinations
of determinants all terms with $x_i\notin\{x_0,x_2,...,x_n\}$ cancel
out and $C_0^{\frac{n}{2}},C_1^{\frac{n}{2}}$ are polynomials in
$x_0,x_2,...,x_n$ of degree $\frac{n}{2}$.

{\em Let $n$ be odd.} We define elements $y_i\in Frac({\cal F}_{od})$ for $i\in\Z$ by 
$$y_{2i}=(f+c)^i,~y_{2i+3}=(f+c)^ig.$$
It is clear that $y_i\in Frac({\cal F}_{od})$ if $i=1,-1,-2,...$ and  $y_0,y_2,...,y_n\in{\cal F}_n$. Moreover, $y_0,y_2,...,y_n$ can be written as linear combinations of 
$x_0,x_2,...,x_n\in{\cal F}_n$.
If $\frac{n+1}{2}$ is even  
we set\footnote{These formulas work for $n>3$. If $n=3$ we set \begin{eqnarray*} 
\tilde{C}_0^{2}&\!\!=\!\!&\det\left((y_0,y_2)^t(y_{-2},y_0)\right),\\ 
\tilde{C}_1^{2}&\!\!=\!\!&\det\left((y_0,y_3)^t(y_0,y_3)\right)-\det\left((y_0,y_2)^t(y_0,a_3y_{2}+b_2y_{3})\right)+\det\left((Q(-c)y_0,y_3)^t(y_{-2},y_0)\right)
\end{eqnarray*}}  
\begin{eqnarray*} 
\tilde{C}_0^{\frac{n+1}{2}}&\!\!=\!\!&\det\left((y_0,y_2,y_4,..,y_{\frac{n+3}{2}})^t(y_{-2},y_0,y_2,..,y_{\frac{n-1}{2}})\right),\\ 
\tilde{C}_1^{\frac{n+1}{2}}&\!\!=\!\!&\ \ \ \det\left((y_0,y_2,..,y_{\frac{n-1}{2}},y_{\frac{n+3}{2}})^t(y_0,y_2,..,y_{\frac{n-1}{2}},y_{\frac{n+3}{2}})\right)\\
&&-\det\left((y_0,y_2,..y_{\frac{n+1}{2}})^t(y_0,y_2,..,y_{\frac{n-1}{2}},a_3y_{\frac{n+1}{2}}+b_2y_{\frac{n+3}{2}})\right)\\
&&-\det\left((Q(-c)y_0,y_2,..,y_{\frac{n-1}{2}},y_{\frac{n+3}{2}})^t(y_{-2},y_0,y_2,y_4,..,y_{\frac{n-1}{2}},y_{\frac{n+3}{2}})\right) \\
&&+\det \left( (Q(-c)y_0,y_2,..,y_{\frac{n+1}{2}})^t (y_{-2},y_0,y_2,y_4,..,y_{\frac{n-1}{2}},a_3y_{\frac{n+1}{2}}+b_2y_{\frac{n+3}{2}}) \right)
\end{eqnarray*}  
and if $\frac{n+1}{2}$ is odd
\begin{eqnarray*}
\tilde{C}_0^{\frac{n+1}{2}}&\!\!\!=\!\!\!&\det\left((y_0,y_2,y_4,..,y_{\frac{n+1}{2}},y_{\frac{n+5}{2}})^t(y_{-2},y_0,y_2,..,y_{\frac{n-3}{2}},y_{\frac{n+1}{2}})\right)\\ &&
-\det\left((y_0,y_2,y_4,..,y_{\frac{n+3}{2}})^t(y_{-2},y_0,y_2,..,y_{\frac{n-3}{2}},b_2y_{\frac{n+1}{2}}+a_3y_{\frac{n-1}{2}})\right),\\ 
\tilde{C}_1^{\frac{n+1}{2}}&\!\!\!=\!\!\!&\det\left((y_0,y_2,..,y_{\frac{n+1}{2}})^t(y_0,y_2,..,y_{\frac{n+1}{2}})\right)
-\det\left((y_{-2},y_0,y_2,y_4,..,y_{\frac{n+1}{2}})^t(Q(-c)y_0,y_2,..,y_{\frac{n+1}{2}})\right) .
\end{eqnarray*}  
In these formulas we use the product in the
algebra $Frac({\cal F}_{od})$ for computing 
products of vector components. 
Therefore
entries of our matrices of the form $v^tw$ are linear combinations of
$y_i\in Frac({\cal F}_{od})$, $i\in\Z$. On the other hand, 
for 
computing
determinants we use the product in $S^{*}(Frac({\cal F}_{od}))$. So
these determinants are polynomials of degree $\frac{n+1}{2}$ in
$y_i,~i\in\Z$. Moreover, it turns out that in our linear combinations
of determinants all terms with 
$y_i\notin\{y_{-2},y_0,y_2,...,y_n\}$
cancel out and $\tilde{C}_0^{\frac{n+1}{2}},\tilde{C}_1^{\frac{n+1}{2}}$
are polynomials in $y_{-2},y_0,y_2,...,y_n$ of degree $\frac{n+1}{2}$.
These polynomials are linear in $y_{-2}$ and therefore can be written
as $\tilde{C}_i^{\frac{n+1}{2}}=A_i+B_iy_{-2}$, $i=0,1$ where $A_i,B_i$
are polynomials in $y_0,y_2,...,y_{n}$.
We set $C^n=A_0B_1-A_1B_0$.

\section{Conclusion} \label{conclude}

In this paper we have constructed nine pairwise compatible quadratic
Poisson structures on a linear space of arbitrary dimension. It seems
that this family of Poisson structures is maximal if the dimension of
linear space is larger than four. We think that the following problems
deserve further investigation:

\begin{itemize}
\item Study the differential and algebraic geometry of these
  compatible Poisson structures to explain geometrically why these
  structures exist and why the number of them is exactly nine.
\item Do there exist other Poisson structures compatible with the one
  in $q_n(\tau)$ where $n\geq 5$?
\item There exist other elliptic Poisson algebras, for example
  $q_{n,k}(\tau)$ where $1\leq k<n$ and $n,~k$ are coprime
  \cite{odes}.
  Functional constructions of these Poisson algebras can be
  found in \cite{odfeig}.  Do there exist Poisson structures
  compatible with the one in $q_{n,k}(\tau)$? Note that
  $q_{n,1}(\tau)=q_n(\tau),~q_{n,n-1}(\tau)$ is trivial so the first
  nontrivial example other than $q_n(\tau)$ is $q_{5,2}(\tau)$.
\end{itemize}

We plan to address these problems elsewhere.

\section*{Appendix: Elliptic and $\theta$-functions} \label{app}
\addcontentsline{toc}{section}{Appendix: Elliptic and $\theta$-functions}

Fix $\tau\in\C$ such that $Im~\tau>0$. Let
$\Gamma=\{k+l\tau;~k,l\in\Z\}\subset\C$ be an integral lattice generated
by 1 and $\tau$. The Weierstrass zeta function is defined as follows:
$$\zeta(z)=\frac{1}{z}+
  \sum_{\omega\in\Gamma\setminus\{0\}}\left(\frac{1}{z-\omega}
                                   +\frac{1}{\omega}
                                   +\frac{z}{\omega^2}\right)$$
The function $\zeta(z)$ is not elliptic but one has
$\zeta(z+\omega)=\zeta(z)+\eta(\omega)$ where $\eta:\Gamma\to\C$ is a
$\Z$-linear function. The functions
$\zeta(z_1-z_2)-\zeta(z_1)+\zeta(z_2)$ and $\zeta^{\prime}(z)$ are
elliptic. Moreover, a function $c_1\zeta(z-u_1)+...+c_m\zeta(z-u_m)$
is elliptic in $z$ if $c_1+...+c_m=0$.

Let $n\in\N$. We denote by $\Theta_n(\tau)$ the space of the entire
functions of one variable satisfying the following relations:
$$f(z+1)=f(z),~~~f(z+\tau)=(-1)^ne^{-2\pi inz}f(z)$$ It is known
\cite{ell} that $\dim\Theta_n(\tau)=n$, every function
$f\in\Theta_n(\tau)$ has exactly $n$ zeros modulo $\Gamma$ (counted
according to their multiplicities), and the sum of these zeros modulo
$\Gamma$ is equal to zero. Let
$\theta(z)=\sum_{\alpha\in\Z}(-1)^{\alpha}e^{2\pi i(\alpha
  z+\frac{\alpha(\alpha-1)}{2}\tau)}$. It is clear that
$\theta(z)\in\Theta_n(\tau)$.  We have $\theta(0)=0$ and this is the
only zero modulo $\Gamma$. Moreover, there exist functions
$\{\theta_{\alpha}(z); \alpha\in\Z/n\Z\}\subset\Theta_n(\tau)$. These
functions are uniquely defined (up to multiplication by a common
constant) by the following identities:
$$ \theta_{\alpha}\left(z+\frac{1}{n}    \right)= e^{ 2\pi i\frac{\alpha}{n}}\theta_{\alpha}(z),
~~~\theta_{\alpha}\left(z+\frac{1}{n}\tau\right)=-e^{-2\pi 
   i(z+\frac{1}{n}-\frac{n-1}{2n}\tau)}\theta_{\alpha}(z)$$ and form a
basis of the linear space $\Theta_n(\tau)$.

Note that $\zeta(z)$ as well as $\theta(z),~\theta_{\alpha}(z)$ are
functions in two variables: $z\in\C$ and modular parameter
$\tau$. Therefore, they can be written as $\zeta(z,\tau)$ and
$\theta(z,\tau),~\theta_{\alpha}(z,\tau)$.


\addcontentsline{toc}{section}{References}

\end{document}